\documentclass[12pt, letterpaper, twoside]{article}
\usepackage[utf8]{inputenc}
\usepackage{amssymb}
\usepackage{amsmath}
\usepackage{amsthm}
\usepackage{comment,xcolor,bm, mathabx}
\newtheorem{lemma}{Lemma}
\newtheorem{definition}{Definition}
\newtheorem{remark}{Remark}
\newtheorem{corollary}{Corollary}

\newtheorem{proposition}{Proposition}
\newtheorem{theorem}{Theorem}
\newtheorem{conjecture}{Conjecture}
\newtheorem{claim}{Claim}
\newcommand{\etal}{\textit{et al.}}
\usepackage[ruled,vlined]{algorithm2e}
\usepackage{setspace}
\newcommand{\C}{{\bf C}}
\newcommand{\HH}{{\bf H}}
\newcommand\ellmodfive{$\ell \equiv 0 \,\, (\mod 5)$}

\newcommand{\A}{{\boldsymbol{A}}}
\newcommand{\bv}{{\boldsymbol{b}}}

\newcommand{\bu}{{\boldsymbol{a}}}

\newcommand{\bq}{{\boldsymbol{q}}}
\newcommand{\1}{{\bf 1}}
\newcommand{\PAF}{\operatorname{PAF}}
\newcommand{\PSD}{\operatorname{PSD}}
\newcommand{\DFT}{\operatorname{DFT}}
\renewcommand{\mod}{\mathrm{mod}\;}

\newcommand{\Z}{\mathbb{Z}}
\newbox\keywbox
\setbox\keywbox=\hbox{\bfseries Keywords:}%
\newcommand\keywords{%
\textbf{Keywords:}\ }
\newcommand{\Ztimes}{\mathbb{Z}_\ell \rtimes \mathbb{Z}^{\times}_\ell}
\newcommand{\ZZtimes}{(\mathbb{Z}_\ell \times \mathbb{Z}_\ell) \rtimes \mathbb{Z}^{\times}_\ell}
\newcommand{\zl}{\mathbb{Z}_\ell}
\newcommand{\zlx}{\zl^\times}
\newcommand{\sdp}{\zl \rtimes \zlx}

\newcommand{\ba}{\mathbf{a}}
\newcommand{\bb}{\mathbf{b}}

\title{Legendre pairs of lengths \ellmodfive} 
\author{
Ilias Kotsireas\thanks{CARGO Lab, Wilfrid Laurier University, Waterloo, Ontario N2L 3C5, Canada}  \,\,\,
Christoph Koutschan\thanks{Supported by the Austrian Science Fund (FWF): F5011-N15. Johann Radon Institute for Computational and Applied Mathematics (RICAM)
Austrian Academy of Sciences, Linz A-4040, Austria}
\,\,\,
Dursun Bulutoglu\thanks{Air Force Institute of Technology, Wright-Patterson Air Force Base, Ohio 45433, USA 
}
\\
David Arquette\footnotemark[3]
\,\,\,
Jonathan Turner\footnotemark[3]
\,\,\,
Kenneth Ryan\thanks{West Virginia University Morgantown, West Virginia 26506, USA}
}
\date{\today}

\begin{document}

\maketitle

\begin{abstract}
\noindent By assuming a type of balance for length $\ell=87$
 and  non-trivial subgroups of multiplier groups of Legendre pairs (LPs) for length $\ell=85$, we find
 LPs
 of these lengths.
 We then study the power spectral density (PSD) values of $m$ compressions of LPs of length $5 m$.
 We also formulate a conjecture for LPs of lengths $\ell \equiv 0$  (mod 5) and demonstrate how it can be used to decrease the search space and storage requirements for finding such LPs.
 The newly found LPs decrease the number of integers in the range $\leq 200$ for which the existence question of LPs remains unsolved from $12$ to $10$.
\end{abstract}
\keywords compressed vector; discrete Fourier transform; Hadamard matrix; periodic autocorrelation function; power spectral density; multiplier group
\section{Introduction}
The \textit{periodic autocorrelation function} (PAF) of a (row) vector \(\ba\in \mathbb{C}^{\ell}\) indexed by \(\zl\) is defined as \(\PAF_\ba (j) = \sum_{i = 0}^{\ell - 1} a_i \overline{a}_{i - j}\), where $\overline{a}_i$ is the complex conjugate of $a_i$.
By using the PAF, we define the concept of a Legendre pair (LP) studied in~\cite{FGS:2001}. 
\begin{definition}\label{def:LP_PAF}
Let $\ba$ and $\bb$ be $\{-1,1\}$ row vectors indexed by $\mathbb{Z}_\ell$.  Then, $(\ba, \bb)$ is an LP if 
\begin{align}
\PAF_{\ba}(j)+\PAF_{\bb}(j)=-2&
\quad \forall j \in \mathbb{Z}_\ell - \{0\}, \label{eqn:pafpaf} \\
\sum_{i=0}^{\ell-1} a_i = \sum_{i=0}^{\ell-1} b_i.& \nonumber
\end{align}
\end{definition}
An LP $(\ba,\bb)$  must satisfy  $\sum_{i=0}^{\ell-1} a_i = \sum_{i=0}^{\ell-1} b_i =1$ or $\sum_{i=0}^{\ell-1} a_i = \sum_{i=0}^{\ell-1} b_i =-1$; see~\cite{Hollon}.


A \textit{circulant shift} of a vector $\ba \in \mathbb{C}^{\ell}$ by $j\in \mathbb{Z}_\ell$, denoted $c_j(\ba)$, is a transformation such that $\left(c_j(\ba)\right)_{i}=a_{{i-j}}$ for each $i\in \mathbb{Z}_\ell$. 
 Let  $\C_{\ba}$ be the  {\em circulant matrix} obtained from a row vector, $\ba$, where the $(j+1)$th row of $\C_{\ba}$ is 
$c_j(\ba)$.

Let $\ba,\bb \in \{-1,1\}^{\ell}$ be an LP such that $\sum_{i=0}^{\ell-1} a_i = \sum_{i=0}^{\ell-1} b_i = 1$, and $\1$ be the length $\ell$ row vector of all $1$s. For a complex number $z=x+iy$, let $\bar{z}=x-iy$ be the complex conjugate of $z$. Then 
\begin{equation}\label{eqn:Hadamard}
\HH=\begin{bmatrix}
-1\phantom{^{\top}} & 
-1\phantom{^{\top}} & 
\phantom{-}\1\phantom{^{\top}} & 
\phantom{-}\1\phantom{^{\top}} \\
-1\phantom{^{\top}} & 
\phantom{-}1\phantom{^{\top}} & 
\phantom{-}\1\phantom{^{\top}} & 
-\1\phantom{^{\top}} \\
\phantom{-}\1^{\top} & 
\phantom{-}\1^{\top} & 
\phantom{-}\C_{\bb}\phantom{^{\top}} & 
\phantom{-}\C_{\ba}\phantom{^{\top}} \\
\phantom{-}\1^{\top} & 
-\1^{\top} & 
\phantom{-}\overline{\C_{\ba}}^{\top} & 
-\overline{\C_{\bb}}^{\top} \\
\end{bmatrix}
\end{equation}
is a $(2\ell+2)\times(2\ell+2)$ {\em Hadamard matrix}, a matrix with entries in $\{\pm1\}$ with orthogonal columns, where for a matrix $\A=(a_{ij})$, $\overline{\A}=(\overline{a_{ij}})$, which will soon become relevant.   It is well known that a Hadamard matrix of order $n$ does not exist if $n$ is not divisible by $4$ and $n>2$. Then $4 \mid (2\ell+2)$, implying that $\ell$ must be odd for an LP of length $\ell$ to exist.  It is conjectured that an LP exists for all odd $\ell$ or equivalently  there is an $n \times n$ Hadamard matrix constructed via LPs whenever $n$ is a multiple of $4$~\cite{Hollon}.

An $n \times n$ matrix with entries from $\{-1,1,-i,i\}$ with orthogonal columns with respect to the complex dot product is called a {\em quaternary (complex) Hadamard matrix}. The order of a quaternary Hadamard matrix must be even~\cite{QuaternaryLP}. 
For even $\ell$, if $\ba,\bb\in \{-1,1,-i,i\}^{\ell}$ satisfy equation~(\ref{eqn:pafpaf}) in  Definition \ref{def:LP_PAF}  and equations
$$ 
\sum_{i=0}^{\ell-1}a_i=1+i,\quad \sum_{i=0}^{\ell-1}b_i=0,
$$ then $\HH$ in equation~(\ref{eqn:Hadamard}) with the upper left $2 \times 2$ submatrix replaced by
$$  
\begin{bmatrix}
-1 & i  \\
-i & 1  \\
\end{bmatrix}
$$
is a quaternary Hadamard matrix, and such $(\ba,\bb)$ is called a {\em quaternary Legendre pair}~\cite{QuaternaryLP}. It is conjectued that a  quaternary LP exists for every even 
$\ell$~\cite{QuaternaryLP}. Since if a quaternary Hadamard matrix of order $n$ exists, then a Hadamard matrix of order $2n$ must also exist~\cite{QuaternaryLP}, proving this conjecture would also prove the {\em Hadamard conjecture}, i.e., a Hadamard matrix of order $n$ exists for every $n$ divisible by $4$.
  

For two groups $N$ and $H$, and an action of $H$ on $N$, let $N\rtimes H$ be the semidirect product of $N$ and $H$ as defined in Rotman~\cite{rotman2010advanced}.   Let $\mathbb{Z}^{\times}_\ell$ be the group of units in $\mathbb{Z}_{\ell}$. Let   $(j,k)i=ki+j$ for 
$(j,k)\in \Ztimes$ and $i \in \mathbb{Z}_{\ell}$. 
Then the group $\Ztimes$ acts on each vector $\ba$ in $\mathbb{C^{\ell}}$ via $(j,k)a_i=a_{(j,k)^{-1}i}$. 
The action of $(0,k):=d_k$ for $k \in \mathbb{Z}^{\times}_\ell$ on $\ba$ is called a {\em decimation}.
 Decimations and cyclic shifts do not commute. 
 In fact, $d_rc_i(d_r)^{-1}=d_rc_id_{r^{-1}}=c_{ir}$. Equivalently, $$c_id_r=d_rc_{ir^{-1}}
 \text{ or } d_rc_k=c_{rk}d_r.  $$

The action of the group $\Ztimes$ on $\ba \in \mathbb{C}^{\ell}$ is used  to define an action of $\ZZtimes$ on LPs
(quaternary LPs) $(\ba,\bb)$ by 
 $((j_1,j_2),k)(\ba,\bb)= ((j_1,k)\ba,(j_2,k)\bb$), 
 for $(j_1,k)\in \Ztimes$ and $(j_2,k)\in \Ztimes$, see~\cite{Hollon,QuaternaryLP}. 
Two pairs of $\{-1,1\}$ or $\{0,1\}$ vectors $(\ba,\bb), (\ba',\bb')$ are \textit{equivalent} if $(\ba,\bb)$ is in the same orbit as $(\ba',\bb')$ or $(\bb',\ba')$  under the action of $\ZZtimes$. Two  pairs of $\{-1,1,-i,i\}$ vectors $(\ba,\bb), (\ba',\bb')$ are \textit{equivalent} if $(\ba,\bb)$ is in the same orbit as $(\ba',\bb')$  under the action of $\ZZtimes$. 
Clearly, if  
$(\ba,\bb)$ is an LP (quaternary LP), and $(\ba',\bb')$ is equivalent to $(\ba,\bb)$, then 
$(\ba',\bb')$ is also  an LP (quaternary 
LP).
The following definition characterizes the inherent symmetries of a vector $\ba \in \mathbb{C}^{\ell}.$ 
\begin{definition}
 For a vector $\ba \in \mathbb{C}^{\ell}$, the group     $$G_{\ba}=\{j \in  \mathbb{Z}_{\ell}^{\times} \mid  (i,j)(\ba) = \ba \text{ for some $i \in \mathbb{Z}_\ell$}\}$$ is called the multiplier group of
$\ba$, and each element of $G_{\ba}$ is called a multiplier of $\ba$. 
\end{definition}
 For \(g \in  \sdp \text{ written as a product of cyclic shifts and decimations and } \bq \in \mathbb{C}^\ell\), 
 the following theorem shows the connection between $G_{\bq}$ and $G_{g\bq}$.
 \begin{theorem}\label{thm:conjugate}
 Let \(g \in  \sdp \text{ and } \bq \in \mathbb{C}^\ell\). 
 Then $G_{g\bq}=gG_{\bq}g^{-1}$ and  $|G_{\bq}|=|G_{g\bq}|$.
 \end{theorem}
 \begin{proof}
Observe that $ r \in G_{\bq}  \Leftrightarrow d_r(\bq)=c_i(\bq)$ for some $i \in \mathbb{Z}_\ell  \Leftrightarrow d_r(\bq)=d_r(g^{-1}(g\bq))=(d_r(g^{-1}))(g\bq)=(c_ig^{-1})(g\bq)=(g^{-1}c_{i'})(g\bq)$ for some  $i' \in \mathbb{Z}_\ell$, where $c_ig^{-1}=g^{-1}c_{i'}$ for some  $i' \in \mathbb{Z}_\ell$
because  $\zl\trianglelefteq \sdp$.  Hence,  $ r \in G_{\bq}  \Leftrightarrow gd_r(g^{-1}(g\bq))=(gd_rg^{-1})(g\bq)=c_{i'}(g\bq)$ for some  $i' \in \mathbb{Z}_\ell$. 
This implies that $gG_{\bq}g^{-1}= G_{g\bq}$.  Then, $|G_{\bq}|=|G_{g\bq}|$ as the map $x \rightarrow gxg^{-1}$ is  one-to-one. 
 \end{proof}

The following theorem appeared as  Theorem 4 
in~\cite{turner2020counting}.
\begin{theorem}\label{thm:djqp0}
 Let $\ell$ be odd and $\bv \in \{0,1\}^\ell$ such that $\sum_{i=0}^{\ell-1}b_i=\delta$ with $\gcd(\delta,\ell)=1$.  Then there is some $i \in \mathbb{Z}_\ell$ such that $\bv'=c_i(\bv)$, and   $$G_{\bv}\leq \widehat{G}_{\bv'}=\{j \in  \mathbb{Z}_{\ell}^{\times} \mid  d_j(\bv') = \bv'\}.$$
\end{theorem}
The next lemma that follows from Theorem~\ref{thm:conjugate} and Theorem~\ref{thm:djqp0}
shows that in Theorem~\ref{thm:djqp0}, $G_{\bv}$ is equal to $\widehat{G}_{\bv'}$.
\begin{lemma}\label{lem:djqp0}
  Let $\ell$ be odd and $\bv \in \{0,1\}^\ell$, such that $\sum_{i=0}^{\ell-1}b_i=\delta$ with $\gcd(\delta,\ell)=1$.  Then there is some $i \in \mathbb{Z}_\ell$ such that $\bv'=c_i(\bv)$, and   $$G_{\bv}= \widehat{G}_{\bv'}=\{j \in  \mathbb{Z}_{\ell}^{\times} \mid  d_j(\bv') = \bv'\}.$$  
\end{lemma}
The following lemma follows from Lemma~\ref{lem:djqp0}.
\begin{lemma}\label{lem:djqppm1}
   Let $\ell$ be odd and $\bv \in \{-1,1\}^\ell$ such that $\sum_{i=0}^{\ell-1}b_i=\delta'$ with $\gcd((\delta'+\ell)/2,\ell)=1$.  Then there is some $i \in \mathbb{Z}_\ell$ such that $\bv'=c_i(\bv)$, and   $$G_{\bv}= \widehat{G}_{\bv'}=\{j \in  \mathbb{Z}_{\ell}^{\times} \mid  d_j(\bv') = \bv'\}.$$  
\end{lemma}
\begin{proof}
Since $(i,j)(\ba) = \ba$ for some $i \in \zl$ if and only if $(i,j)(\ba+\1)/2 = (\ba+\1)/2$ for some $i \in \zl$. Hence, $$G_{\bv}=G_{\frac{\bv+\1}{2}}.$$ Since, $\sum_{i=0}^{\ell-1}(b_i+1)/2=(\delta'+\ell)/2$
 and  $\gcd((\delta'+\ell)/2,\ell)=1$, by Lemma~\ref{lem:djqppm1}  there is some $i \in \mathbb{Z}_\ell$ such that $\bv'=c_i((\bv+\1)/2)=(c_i(\bv)+\1)/2$, 
 and   $$G_{\bv}=G_{\frac{\bv+\1}{2}}= \widehat{G}_{\bv'}=\widehat{G}_{\frac{c_i(\bv)+\1}{2}}.$$ 
 Now, $$\widehat{G}_{\frac{c_i(\bv)+\1}{2}}=\widehat{G}_{c_i(\bv)}=\widehat{G}_{\bv'} $$ follows from the definition of 
 $\widehat{G}_{\bv'} $.
\end{proof}
The following corollary simplifies the search for an LP $(\bu,\bv)$ with multiplier groups 
$G_{\bu}$ and $G_{\bv}$. 
\begin{corollary}\label{cor:GuhatGu}
An LP $(\bu,\bv)$ with multiplier groups $G_{\bu}$ and $G_{\bv}$ exists if and only if an equivalent LP $(\bu',\bv')$ with multiplier groups $\widehat{G}_{\bu'}=G_{\bu}$ and $\widehat{G}_{\bv'}=G_{\bv}$ exists. Hence, a search for all non-equivalent LPs $(\bu,\bv)$ with multiplier groups $G_{\bu}$ and $G_{\bv}$
can be implemented by searching for all non-equivalent  LPs $(\bu',\bv')$ with multiplier groups $\widehat{G}_{\bu'}$ and $\widehat{G}_{\bv'}$.
\end{corollary}
\begin{proof}
 For an LP $(\bu,\bv)$, $\sum_{i=0}^{\ell-1}a_i=\sum_{i=0}^{\ell-1}b_i=\delta'=\pm1$ and $\gcd((\ell\pm1)/2,\ell)=1$.   Hence, by Lemma~\ref{lem:djqppm1},  there exits $i_1,i_2 \in \zl$ such that $G_u=\widehat{G}_{c_{i_1}(\bu)}$ and $G_v=\widehat{G}_{c_{i_2}(\bv)}$, and $(c_{i_1}(\bu),c_{i_2}(\bv))$ is an LP equivalent to $(\bu,\bv)$ with multiplier groups $\widehat{G}_{c_{i_2}(\bu)}=G_{\bu}$, and   $\widehat{G}_{c_{i_2}(\bv)}=G_{\bv}$. 
\end{proof}
Our repository of LPs, provides many examples of LPs  $(\bu,\bv)$, with $G_u \neq G_v$.
However, for LPs  $(\bu,\bv)$ in our repository with non-trivial multiplier groups, more often than not $G_u=G_v$.
On the other hand,  for most LPs in our repository, $G_\bu=G_\bv = \{ 1 \} $. 

Throughout the rest of the paper, let the row vectors $\ba, \bb\in \mathbb{R}^{\ell}$, $\ell$ be odd,  and  \(w = e^{2 \pi i / \ell}\). The \textit{discrete Fourier transform} of $\ba$ is  $\DFT_\ba (j)   := \sum_{r = 0}^{\ell - 1} w^{j r} a_r$  and 
the \textit{power spectral density} of $\ba$ is  \(\PSD_{\ba}(j):=|\DFT_\ba (j)|^2 \) for $j \in \mathbb{Z}_{\ell}$.
It is shown in~\cite{FGS:2001} that a pair of row vectors \(\ba,\bb\in \{-1,1\}^\ell\) form an LP of length~$\ell$ if and only if
\begin{gather}
\sum_{i=0}^{\ell-1} a_i = \sum_{i=0}^{\ell-1} b_i =\pm 1, \label{eqn:LPpm}\\
\PSD_{\ba}(j) + \PSD_{\bb}(j) = 2\ell+2 \quad \forall  j \in \zl - \{0\}.\label{eqn:LPpm2}
\end{gather}
If \(\ba \in \{-1,1\}^\ell\), then \((\ba+\1)/2 \in \{0,1\}^\ell\). It is plain to see that
\begin{equation}\label{eqn:PSDone}
\PSD_{\frac{\ba+\1}{2}}(j)=\frac{\PSD_{\ba}(j)}{4} \quad \forall  j \in \zl - \{0\},
\end{equation}
and
\begin{equation}\label{eqn:PSDone2}
\PSD_{\frac{\ba+\1}{2}}(0)=|\DFT_\frac{\ba+\1}{2} (0)|^2 =
\left(\frac{\sum_{i=0}^{\ell-1} a_i+\ell}{2}\right)^2=\left(\frac{\ell\pm1}{2}\right)^2.
\end{equation}
The following lemma follows from equations~(\ref{eqn:LPpm}),~(\ref{eqn:LPpm2}),~(\ref{eqn:PSDone}), and~(\ref{eqn:PSDone2}).
\begin{lemma}\label{lem:LPone}
Let \(\ba,\bb\in \{-1,1\}^\ell\) form an LP of length~$\ell$.    
  Then 
\begin{gather*}
\sum_{i=0}^{\ell-1} \frac{a_i+1}{2} = \sum_{i=0}^{\ell-1} \frac{b_i+1}{2} =\frac{\ell \pm 1}{2}, \label{eqn:LPone}\\
\PSD_{\frac{\ba+\1}{2}}(j) + \PSD_{\frac{\bb+\1}{2}}(j) = \frac{\ell+1}{2} \quad \forall  j \in \zl - \{0\}.\label{eqn:LPone2}
\end{gather*}
\end{lemma}
The following theorem from~\cite{FGS:2001} determines the relation between $\PSD_{\ba}(j)$ and $\PAF_{\ba}(j)$.
\begin{theorem}\label{thm:Wiener}
[Wiener-Khinchin Theorem] The PSD of  $\ba \in \mathbb{R}^{\ell}$ is equal to the
DFT of its PAF, i.e., 
\begin{equation*}
\PSD_{\ba}(k) = \sum_{j=0}^{\ell-1} \PAF_{\ba}(j)w^{jk}\quad \text{for } k \in \mathbb{Z}_{\ell}.
\end{equation*}
Moreover, the PAF of $\ba$ is equal to the inverse DFT of $\ba$'s PSD, i.e., 
\begin{equation*}
\PAF_{\ba}(j) = \frac{1}{\ell}\sum_{k=0}^{\ell-1} \PSD_{\ba}(k)w^{-jk}\quad \text{for } j \in \mathbb{Z}_{\ell}.
\end{equation*}
\end{theorem}
The following corollary follows directly from Theorem~\ref{thm:Wiener} and equation~(\ref{eqn:PSDone2}).
\begin{corollary}\label{cor:sumpaf}
Let \(\ba, \bb \in\{-1,1\}^\ell\) form an LP of length~$\ell$. Then
$$1=\PSD_{\ba}(0) = \sum_{j=0}^{\ell-1} \PAF_{\ba}(j),$$ and
$$\left(\frac{\ell\pm1}{2}\right)^2=\PSD_{\frac{\ba+\1}{2}}(0) = \sum_{j=0}^{\ell-1} \PAF_{\frac{\ba+\1}{2}}(j).$$
\end{corollary}

Throughout the paper, WLOG, we assume that $\sum_{i=0}^{\ell-1} a_i = \sum_{i=0}^{\ell-1} b_i =1$.
Also, let 
$e_2(\ba)=\sum_{i<j}a_ia_j$ be the {\it elementary symmetric  polynomial of degree $2$} and $p_2(\ba)=\sum_{i}a_i^2$ be the {\it power sum symmetric polynomial of degree $2$}.

For $\ba=[a_0,\ldots,a_{\ell-1}]$ and $\ell=nm$ for some positive integers $n,m$,  define $ A_j = \sum_{i=0}^{m-1} a_{ni+j}$ for $j=0,\ldots,n-1$. The vector $\mathcal{A}_m = [A_0 , \ldots , A_{n-1} ]$ is called the {\it $m$-compression} of $\ba$~\cite{DK:DCC:2015}. Throughout the paper, let $\mathcal{A}_m = [A_0 , \ldots , A_{n-1} ]$,  $\mathcal{B}_m = [B_0 , \ldots , B_{n-1} ]$
 be the $m$-compressions of $\ba$, $\bb$.
 
The following theorem from~\cite{TKBG:DCC:2021}
shows how the DFT and the PSD of the $m$-compression of a vector $\ba$ are related to the DFT and the PSD of $\ba$.
\begin{theorem}\label{thm:compression}
Let $\ba$ be a vector of length $\ell= nm$, and $\mathcal{A}_m$ be the $m$-compression of $\ba$. Then 
$\DFT_{\mathcal{A}_m} (j)=\DFT_\ba (mj)$ and  $\PSD_{\mathcal{A}_m} (j)=\PSD_\ba (mj)$ for $j \in \{0,\ldots,n-1\}$.
\end{theorem}
Compression of complementary 
vectors has proved to be a valuable tool for finding several
previously unknown complementary vectors (vectors whose  PAF vectors' $j$th entries sum to a constant for $j \neq 0$) in the past decade~\cite{DK:DCC:2015,KK:JCD:2021,QuaternaryLP,TKBG:DCC:2021}.  Compression-based search algorithms are based on a two-step process. 
In the first step, several candidate compression vectors are computed. The second step involves searching for decompressions of the candidate compressed vectors from the first step.


In Section~\ref{sec:mod5}, we find restrictions on the PSD values of LP vectors computed using only the $5$th primitive  roots of unity. We then corroborate our theoretical restrictions by computing all PSD values which are only based on the  $5$th primitive roots of unity for all LP vectors that we possess. Then, based on these restrictions, we 
develop Conjecture~\ref{conj:Ilias} which applies to $\ell\equiv 0$ (mod 5) cases. We then show how this conjecture can be used to prune the search space.
Moreover, we provide evidence for Conjecture~\ref{conj:Ilias} by confirming that it is valid for all $\ell \leq 85$ such that $\ell\equiv 0$ (mod $5$) (with the sole exception of $\ell = 75$). 

In Section~\ref{sec:new85}, we first discuss the limitations of the methods used by~\cite{TKBG:DCC:2021} for finding LPs of length $\ell>77$ and how those limitations can be overcome.
Then, by assuming that the sought after LP $(\bu,\bv)$ satisfies $\{1,69\}\leq \widehat{G}_{\bu}$,  
$\{1,69\}\leq \widehat{G}_{\bv}$, and using  the method in Section 4.1.1 of~\cite{KK:JCD:2021} with $H_1=\{ 1, 69 \}$, we find the first known examples of LPs of length $85$ that satisfy this property.
The length $\ell=85$ has been the smallest open length case.
Additionally,  we show how to exploit a different concept 
of balance applied to LPs of composite order. As an application, we find  the first examples of LPs of length $87$. 

In Section~\ref{sec:search115}, we describe the partial searches that have been implemented for LPs of length $115$.

In Section~\ref{sec:conclusion}, we provide the current list of ten $\ell$ values less than $200$ for which the LP existence problem remains unsolved.
\section{LPs of lengths \ellmodfive}\label{sec:mod5}
For a vector of length $\ell=n m$, the following lemma  expresses a sum of PAFs in terms of the PAF of the compressed vector $\mathcal{A}_m$.
\begin{lemma}\label{lem:Koutschan}
Let $\ba$ be a length $\ell=n m$ vector for positive odd integers $n$ and $m$. Then
\begin{equation*}
  \sum_{j=0}^{m-1} {\rm PAF}_a(n j+k) = {\rm PAF}_{{\cal A}_m}(k).
\end{equation*}
\end{lemma}
\begin{proof}
\begin{align*}
  \sum_{j=0}^{m-1} {\rm PAF}_a(n j+k) &=
  \sum_{j=0}^{m-1} \sum_{i=0}^{\ell-1} a_i a_{i-n j-k} =
  \sum_{j=0}^{m-1} \sum_{i=0}^{m-1} \sum_{r=0}^{n-1} a_{n i+r} a_{n i-n j+r-k} \\ &=
  \sum_{r=0}^{n-1} \sum_{i=0}^{m-1} \sum_{j=0}^{m-1} a_{n i+r} a_{n j+r-k} =
  \sum_{r=0}^{n-1} \left(\sum_{i=0}^{m-1} a_{n i+r}\right)
  \left(\sum_{j=0}^{m-1} a_{n j+r-k}\right) \\ &=
  \sum_{r=0}^{n-1} A_r A_{r-k} = {\rm PAF}_{{\cal A}_m}(k).
  \qedhere
\end{align*}
\end{proof}
For a vector of length $\ell=5m$, the following lemma determines ${\rm PSD}_{\ba}(rm)$ in terms of ${\rm PAF}_{\ba}(j)$ for $j \in \mathbb{Z}_{\ell}$.
\begin{lemma}\label{lem:Dursun}
Let $m \in \mathbb{Z}_{\geq 1}$, $r \in \{1,\ldots,4\}$, and $\ba$ be a vector of length $\ell=5m$. Then 
\[
  {\rm PSD}_{\ba}(rm)
  =\sum_{j=0}^{m-1} {\rm PAF}_{\ba}(5j)
  -\sum_{k=1}^2 \frac{1+(-1)^{k+\left\lfloor \frac{r}{2}\right\rfloor}\sqrt{5}}{2}
  \sum_{j=0}^{m-1}{\rm PAF}_{\ba}(5j+k).
\]
\end{lemma}
\begin{proof}
Let $w=e^{2\pi i/\ell}$, $w'=w^m$. 
By Theorem~\ref{thm:Wiener}, 
\begin{align*}
\begin{split}
{\rm PSD}_{\ba}(rm)=&\sum_{j=0}^{\ell-1}
{\rm PAF}_{\ba}(j)w^{rmj}=\sum_{j=0}^{\ell-1}{\rm PAF}_{\ba}(j)(w')^{rj}=
\sum_{k=0}^4\sum_{j=0}^{m-1}{\rm PAF}_{\ba}(5j+k)(w')^{rk}\\
=&\sum_{j=0}^{m-1}{\rm PAF}_{\ba}(5j)+\sum_{k=1}^2
2\cos \left (\frac{2\pi rk}{5}\right)\sum_{j=0}^{m-1}{\rm PAF}_{\ba}(5j+k)\\
=&\sum_{j=0}^{m-1}
{\rm PAF}_{\ba}(5j)+\sum_{k=1}^2\left(\frac{-1+
(-1)^{\left\lfloor(rk\ {\mod5})/2\right\rfloor}\sqrt{5}}{2}\right)\sum_{j=0}^{m-1}{\rm PAF}_{\ba}(5j+k),\\
\end{split}
\end{align*}
from which the assertion follows directly.
\end{proof}

\begin{corollary}\label{cor:Dursun}
For a pair of vectors $(\ba,\bb)$ of length 
$\ell = 5 m$,  the LP constraints
\[
  {\rm PSD}_{\ba}(rm)+{\rm PSD}_{\bb}(rm)=2\ell+2
\]
are satisfied for $r \in \{1,\ldots,4\}$ if and only if 
\[
  \sum_{j=0}^{m-1} 
  \Bigl({\rm PAF}_{\ba}(5j)+{\rm PAF}_{\bb}(5j)\Bigr)
  -\sum_{j=0}^{m-1}\sum_{k=1}^2\frac{{\rm PAF}_{\ba}(5j+k)+{\rm PAF}_{\bb}(5j+k)}{2}=2\ell+2,
\]
and 
\[
  \sum_{j=0}^{m-1}
  \sum_{k=1}^2(-1)^k\Bigl({\rm PAF}_{\ba}(5j+k)+{\rm PAF}_{\bb}(5j+k)\Bigr)=0.
\]
\end{corollary}

The following corollary follows from Lemmas~\ref{lem:Koutschan} and~\ref{lem:Dursun}.
\begin{corollary}\label{cor:PSDarm}
Let $\ba$ be a vector of length $5m$ for some $m \in \mathbb{Z}_{\geq 1}$. Let the $m$-compression of $\ba$ be ${\cal A}_m = [A_0,A_1,A_2,A_3,A_4]$. 
 Then, 
\begin{align}\label{eqn:Koukouvinos}
\begin{split}
     \PSD_{\ba}(rm) &=  p_2({\cal A}_m) -\frac{1}{2}e_2({\cal A}_m) +(-1)^{\left\lfloor \frac{r}{2}\right\rfloor} \sqrt{5} \left(\frac{ \PAF_{{\cal A}_m}(1) - \PAF_{{\cal A}_m}(2) }{2}\right) \\
\end{split}
\end{align}
for $r=1,\ldots,4$.
\end{corollary}

\begin{proposition}\label{prop:n1n2}
Let  $(\ba,\bb)$ be an LP of length $\ell = 5 m$ with $m$ odd. Then there exist integers  $n_1,n_2\in \mathbb{Z}_{\geq 0}$ with $n_1+n_2=2\ell+2$ and $x\in\Z$ such that
\begin{align}\label{eqn:sqrtx}
\begin{split}
    \PSD_{\ba}(rm) &= n_1 + \displaystyle (-1)^{\left\lfloor \frac{r}{2}\right\rfloor} \sqrt{5}  x, \\
    \PSD_{\bb}(rm) &= n_2 - \displaystyle (-1)^{\left\lfloor \frac{r}{2}\right\rfloor} \sqrt{5}  x
    \end{split}
\end{align}
for $r=1,\ldots,4,$
where 
\begin{equation}\label{eqn:PAFab}
x= \frac{ \PAF_{{\cal A}_m}(1) - \PAF_{{\cal A}_m}(2) }{2}= -\frac{ \PAF_{{\cal B}_m}(1) - \PAF_{{\cal B}_m}(2) }{2}.
\end{equation}

\end{proposition}
\begin{proof}
Equations~(\ref{eqn:sqrtx}) directly follow from equation (\ref{eqn:Koukouvinos}) and the fact that an LP  $(\ba,\bb)$ of length $\ell = 5 m$ satisfies ${\rm PSD}_{\ba}(rm)+{\rm PSD}_{\bb}(rm)=2\ell+2$ for $r=1,\ldots,4$.
Since by our general assumption the vectors $\ba,\bb$ satisfy
$\sum_{i=0}^{\ell-1} a_i = \sum_{i=0}^{\ell-1} b_i =1$,
the same holds true for their $m$-compressions, that is
$\sum_{i=0}^4A_i=\sum_{i=0}^4B_i=1$.
Consequently, $\left(\sum_{i=0}^4A_i\right)^2=p_2({\cal A}_m)+2e_2({\cal A}_m)=1$,
and analogously for ${\cal B}_m$. Then by Corollary~\ref{cor:PSDarm},
\begin{equation}\label{eqn:1a2}
\begin{alignedat}{2}
   n_1 &= p_2({\cal A}_m) -\frac{1}{2}e_2({\cal A}_m) &&= 
   \frac{1}{4} \Bigl( 5 p_2({\cal A}_m) -1\Bigr), \\
   n_2 &= p_2({\cal B}_m) -\frac{1}{2}e_2({\cal B}_m) &&= 
   \frac{1}{4} \Bigl( 5 p_2({\cal B}_m) -1\Bigr).
\end{alignedat}
\end{equation}
This shows that $n_1$ and $n_2$ are integers, because $p_2({\cal A}_m)$ is the
sum of five odd squares and hence equals $1\ {\mod4}$.
To see that $x$ is an integer, note that an odd~$m$ implies that $\PAF_{{\cal A}_m}(k)$ and $\PAF_{{\cal B}_m}(k)$ for $k=1,2$ are all odd. Equations~(\ref{eqn:PAFab}) follow from Corollary~\ref{cor:PSDarm} and equations~(\ref{eqn:sqrtx}).

We will next  prove that $0\leq n_1$. By symmetry, this will also give us $0\leq n_2$. 
By equations~(\ref{eqn:Koukouvinos}) and~(\ref{eqn:1a2})
$$\frac{4n_1}{25}=\frac{p_2(\mathcal{A}_m)}{5}-\frac{1}{25}= \frac{p_2(\mathcal{A}_m)}{5}-\left(\frac{\sum_{i=0}^4A_i}{5}\right)^2.$$
Hence, $4n_1/25$ is the population variance of the numbers in $\{A_0,A_1,A_2,A_3,A_4\}$ implying that $0\leq n_1$.
\end{proof}
Experimental evidence gathered for $\ell = 5,15,25,35,45,55,65,85$ indicates  there are LPs of these orders such that 
\begin{align}\label{eqn:5p2}
\begin{split}
  n_1= \frac{1}{4} \Bigl( 5 p_2({\cal A}_m) -1\Bigr)& = \ell + 1=5m+1, \\
  n_2= \frac{1}{4} \Bigl( 5 p_2({\cal B}_m) -1\Bigr)& = \ell + 1=5m+1. 
\end{split}
\end{align}
Solving equations~(\ref{eqn:5p2}) for $p_2({\cal A}_m)$ and $p_2({\cal B}_m)$, we get that equations~(\ref{eqn:5p2}) are equivalent to 
\begin{equation*}\label{eqn:p2a2}
  p_2({\cal A}_m) = p_2({\cal B}_m) = 4m + 1.
\end{equation*}

The next two lemmas determine restrictions on $n_1,n_2$, and $x$ in Proposition~\ref{prop:n1n2}.
\begin{lemma}\label{lem:n1x}
Each of $n_1,n_2$, and $x$ in Proposition~\ref{prop:n1n2} must be even.
\end{lemma}
\begin{proof}
By symmetry, it suffices to prove the result for $n_1$ and $x$. Let $\ba$ be as in Proposition~\ref{prop:n1n2}. By replacing $\ba$ with $(\ba+\1)/2$ in Lemma~\ref{lem:Dursun}, and by equation~(\ref{eqn:PSDone}) we get that 
\[
  {\rm PSD}_{\frac{\ba+\1}{2}}(m) = \frac{\alpha_1}{2}+\frac{\alpha_2}{2}\sqrt{5}=\frac{n_1+\sqrt{5}x}{4},
\]
for some $\alpha_1,\alpha_2 \in \mathbb{Z}$. Hence, by the linear independence of $1$ and $\sqrt{5}$ in the field extension $\mathbb{Q}[\sqrt{5}]$, we have $\alpha_1=n_1/2$ and $\alpha_2=x/2$, implying that both $n_1$ and $x$ must be even.
\end{proof}

\begin{lemma}\label{lem:n1px}
Let $n_1,n_2$, and $x$ be as in Proposition~\ref{prop:n1n2}. Then $n_i+x\equiv 0 \,\, (\mod 4)$ for $i=1,2$.
\end{lemma}
\begin{proof}
By symmetry, it suffices to prove the result for $n_1$ and $x$. Let $w'=e^{2\pi i/5}$  
and  $\ba$ be as in Proposition~\ref{prop:n1n2}. Then, by the proof of Lemma~\ref{lem:Dursun},
\begin{align*}
\begin{split}
{\rm PSD}_{\ba}(m)=&
\sum_{k=0}^4\sum_{j=0}^{m-1}{\rm PAF}_{\ba}(5j+k)(w')^k=n_1+\sqrt{5}x.
\end{split}
\end{align*}
Now, by equation~(\ref{eqn:PSDone}),
\begin{align*}
\begin{split}
{\rm PSD}_{\frac{\ba+\1}{2}}(m)=&\sum_{k=0}^4\sum_{j=0}^{m-1}{\rm PAF}_{\frac{\ba+\1}{2}}(5j+k)(w')^k=\frac{n_1+\sqrt{5}x}{4}.
\end{split}
\end{align*}
 Since  ${\rm PAF}_{(\ba+\1)/2}(5j+k)={\rm PAF}_{(\ba+\1)/2}(-5j-k)$ for $k \in \{1,2\}$,
\begin{align*}
{\rm PSD}_{\frac{\ba+\1}{2}}(m)= & S_0+S_1\left(w'+\overline{w'}\right)+S_2\left((w')^2+\left(\overline{w'}\right)^2\right), 
\end{align*}
where $S_0=\sum_{j=0}^{m-1}{\rm PAF}_{(\ba+\1)/2}(5j)$ and $S_i=\sum_{j=0}^{m-1}{\rm PAF}_{(\ba+\1)/2}(5j+i)$ for $i=1,2$. 
Observe that $S_0,S_1,S_2 \in \mathbb{Z}_{\geq 1}$.
Then, since $(w')^2+\left(\overline{w'}\right)^2=\left(w'+\overline{w'}\right)^2-2\in \mathbb{Z}\left[w'+\overline{w'}\right]$,
\begin{align*}
{\rm PSD}_{\frac{\ba+\1}{2}}(m)
=S_0+S_1\left(w'+\overline{w'}\right)+S_2\left((w')^2+\left(\overline{w'}\right)^2\right)\in \mathbb{Z}\left[w'+\overline{w'}\right].
\end{align*}
Then since $\sqrt{5}=2w'+2\overline{w'}+1$, 
\begin{align*}
{\rm PSD}_{\frac{\ba+\1}{2}}(m)
=\frac{n_1+\sqrt{5}x}{4}=\frac{n_1+x+2x\left(w'+\overline{w'}\right)}{4}\in \mathbb{Z}\left[w'+\overline{w'}\right].
\end{align*}
Hence, $n_1+x \equiv 0 \, \, \, (\mod 4)$.
\end{proof}
The following general lemma  will be used to find another restriction on $n_1$ and $n_2$ in Proposition~\ref{prop:n1n2}. 
\begin{lemma}
\label{lem:Turner}
Let $\ba$ be a vector of length $\ell=mn$ for some positive integers $m$ and $n$. Then $${\sum_{k=0}^{n-1}\PSD_{\ba}(km)=n\sum_{x=0}^{m-1}{\rm PAF}_{\ba}(xn)}=n\, {\rm PAF}_{{\cal A}_m}(0).$$
\end{lemma}
\begin{proof} 
 Since $\PSD_a(km) = \PSD_{\mathcal{A}_m}(k)$, by applying the inverse Fourier transform to
 $[\PSD_{\mathcal{A}_m}(0)\ldots\PSD_{\mathcal{A}_m}(n-1) ]$, we get $[\PAF_{\mathcal{A}_m}
 (0),\ldots,\PAF_{\mathcal{A}_m}(n-1)]$, which at $0$ is what the lemma claims
 (cf. Theorem~\ref{thm:Wiener}).
\end{proof} 
The following lemma  provides an additional restriction on the values of $n_1, n_2$ in Proposition~\ref{prop:n1n2}.
\begin{lemma}\label{lem:n1n2mod6}
Let $\ba,\bb, n_1,$ and $n_2$ be as in Proposition~\ref{prop:n1n2}. Then
 $ n_i\equiv 6 \,\, (\mod 10)$ for $i=1,2$.
\end{lemma}
\begin{proof}
By symmetry, it suffices to prove the result for $i=1$. By Lemma~\ref{lem:Turner} and Proposition~\ref{prop:n1n2} $$\sum_{k=0}^{4}\PSD_{\ba}(km)=1+4n_1=5\left(\sum_{x=0}^{m-1}{\rm PAF}_{\ba}(5x)\right)\equiv 0 \,\, \, (\mod 5).$$
Alternatively, we can employ $p_2({\cal A}_m)=4s+1$ (see the proof of 
Proposition~\ref{prop:n1n2}) to deduce that
$n_1 =  \bigl( 5 p_2({\cal A}_m) -1\bigr)/4=5s+1$ for some integer~$s$.
Now the result follows from the Chinese remainder theorem since $n_1$ must be even by Lemma~\ref{lem:n1x}. 
 \end{proof}
 The following corollary  provides upper and lower bounds on $n_1, n_2$ in Proposition~\ref{prop:n1n2}.
 \begin{corollary}\label{cor:boundn1n2}
Let $\ell=5m$, and let $\ba$, $\bb$, $n_1$, and $n_2$ be as in Proposition~\ref{prop:n1n2}. Then
 $6 \leq n_i \leq 2\ell-4 $ 
  for $i=1,2$.
\end{corollary}
\begin{proof} 
By Proposition~\ref{prop:n1n2}, we get
$0 \leq n_i \leq 2\ell+2=10m+2$ for $i=1,2$. By Lemma~\ref{lem:n1n2mod6}, we get $6 \leq n_i \leq 2\ell-4=10m-4$ for $i=1,2$.
\end{proof}
The lower bound for $n_i$ in Corollary~\ref{cor:boundn1n2} is achieved by taking
${\mathcal A}_{m}=\left[1,1,1,-1,-1\right].$

The following corollary now follows from Proposition~\ref{prop:n1n2}, 
Corollary~\ref{cor:boundn1n2}, and Lemma~\ref{lem:n1px}.
\begin{corollary}\label{cor:k1plusk2}
Let $\ell=5m$, and let $\ba$, $\bb$, $n_1$, $n_2$, and $x$ be as in Proposition~\ref{prop:n1n2}.
Then there exist non-negative integers $k_1,k_2$ such that $n_i=10k_i+6$ for $i=1,2$, $k_1+k_2=m-1$, 
and $2k_i+2\equiv x$  (mod $4$).
\end{corollary}
Corollary~\ref{cor:k1plusk2} summarizes all the constraints obtained in this section regarding $n_1$, $n_2$, and $x$ in equation~(\ref{eqn:sqrtx}). In fact, all the constraints regarding $n_1,n_2$, and $x$ that appear before Corollary~\ref{cor:k1plusk2} are implied by  Corollary~\ref{cor:k1plusk2}.

 Exhaustive searches for $\ell = 5, 15, 25$ have revealed  that most LPs of these lengths have $n_1=n_2=\ell+1$, where $n_1,n_2$ are as in Proposition~\ref{prop:n1n2}. 
 The case $n_1=n_2=\ell+1$ corresponds to 
 taking $k_1=k_2=(m-1)/2$ in Corollary~\ref{cor:k1plusk2}. 
\noindent Non-exhaustive searches for larger odd values of $\ell$ which are multiples of $5$ have revealed the same pattern, i.e.,
the standard relationship $\PSD_{\ba}(m) + \PSD_{\bb}(m) = 2\ell +2$ is materialized (often, but not always) in the same  manner. These observations give rise to the following conjecture.
\begin{conjecture}\label{conj:Ilias}
For every odd positive integer \ellmodfive, there exists at least one LP $(\ba,\bb)$ of length $\ell = 5  m$ such that 
\begin{align}\label{eqn:Ilias}
\begin{split}
\PSD_{\ba}(rm) &=\ell+ 1 + (-1)^{\left\lfloor\frac{r}{2}\right\rfloor} \displaystyle\sqrt{5}x \\
\PSD_{\bb}(rm) &=\ell+ 1 - (-1)^{\left\lfloor\frac{r}{2}\right\rfloor} \displaystyle\sqrt{5}x
\end{split}
\end{align} 
for $r=1,\ldots,4$, and for some nonnegative  integer   $x\equiv -(\ell+1)\equiv \ell+1  \pmod 4$. 
\end{conjecture}
By formula~(\ref{eqn:sqrtx}), Conjecture~\ref{conj:Ilias} holds for one $r$ if and only if $n_1=n_2.$ Again, by formula~(\ref{eqn:sqrtx}), if Conjecture~\ref{conj:Ilias} holds for one $r$, then it also holds for the other $r$s. 
Conjecture~\ref{conj:Ilias} states that for each odd positive \ellmodfive,  there exists a length $\ell$  LP $(\ba,\bb)$ such that $n_1=n_2=\ell+1$ in Proposition~\ref{prop:n1n2}. That is,
the constant $2\ell + 2$ can be
distributed in a ``balanced'' manner between $\PSD_{\ba}(rm)$ and $\PSD_{\bb}(rm)$ for $r \in \{1,\ldots,4\}$. On the other hand,  there exist LPs $(\ba,\bb)$ of length $\ell = 5  m$ which do not satisfy $n_1=n_2$.
Table~\ref{tab:Conjecture} 
provides computational evidence for Conjecture~\ref{conj:Ilias}. Cases $m=1,3,5$ are based on complete---and the remaining cases are based on partial---classifications of all non-equivalent LPs, where equivalent LPs satisfying Conjecture~\ref{conj:Ilias} have the same $x$ value. This is because 
$$
{\rm PSD}_{(j_1,k)\ba}(rm)= {\rm PSD}_{\ba}(krm)
$$
for each $(j_1,k) \in \Ztimes$ and $r \in \{1,\ldots,4\}$.
\begin{table}[ht!]
\centering
    \begin{tabular}{c|c|c}
        $m$ & $\ell  = 5  m$ & $x$ \\
        \hline 
        1 & 5 & $2$ \\
        3 & 15 & $0,4$ \\
        5 & 25 & $2,6,10$ \\
        7 & 35 & $0,8,16$ \\
        9 & 45 & $2,6,10,14,18$ \\
        11 & 55 & $4,8,12,20,24$ \\
        13 & 65 & $14$ \\
        \hline
        17 & 85 & $18$ \\
        \hline 
    \end{tabular}
    \caption{Computationally verified cases for Conjecture~\ref{conj:Ilias} with their corresponding $x$ value(s)}
    \label{tab:Conjecture}
\end{table}

The following proposition provides necessary and sufficient constraints on the $m$-compressed vectors  $(\mathcal{A}_{m},$ $ \mathcal{B}_{m})$ of an LP $(\ba,\bb)$ of length $\ell = 5m$
to satisfy equations~(\ref{eqn:Ilias}). 
\begin{proposition}\label{proposition_1}
The  $m$-compressed vectors  ${\mathcal A}_{m}=[A_0,\ldots,A_4]$ and ${\mathcal B}_{m}=[B_0,\ldots,B_4]$ of an LP $(\ba,\bb)$ of length $\ell = 5 m$ satisfy the conditions in Conjecture~\ref{conj:Ilias}
if and only if 
  \begin{equation}\label{SOS_Diophantine_equations}
  p_2(\mathcal{A}_m)={\rm PAF}_{\mathcal{A}_{m}}(0) = p_2(\mathcal{B}_m)={\rm PAF}_{\mathcal{B}_{m}}(0) = 4m+1.
  \end{equation}
\end{proposition}
\begin{proof}
This result follows from Proposition~\ref{prop:n1n2} and  equation~(\ref{eqn:1a2}). 
\end{proof} 
 Equations~(\ref{SOS_Diophantine_equations}) in Proposition~\ref{proposition_1} drastically reduce the search space for an LP $(\ba,\bb)$ satisfying the conditions in Conjecture~\ref{conj:Ilias}. For such an LP, we are only interested in all-odd solutions to equations~(\ref{SOS_Diophantine_equations})
since $A_i, B_i$ are odd numbers for $i = 0,\ldots,4$. 
Another consequence of Proposition~\ref{proposition_1} is that the alphabet of the possible $m$-compressed vectors $({\mathcal A}_{m}, {\mathcal B}_{m})$ is also significantly truncated. More specifically,
while the full alphabet for candidate $m$-compressed vectors 
is  $$\{ -m, -(m-2), \ldots, -1, +1, \ldots, (m-2), m \},$$
by Proposition~\ref{proposition_1},  
 $A_i$ is an odd integer with $|A_i|\leq \sqrt{4m-3}<2\sqrt{m}$ for $i \in \mathbb{Z}_{\ell}$. 

 Next, we show 
Conjecture~\ref{conj:Ilias} cannot be ruled out by the necessary constraints~(\ref{SOS_Diophantine_equations}) 
and 
$
\sum_{i=0}^4A_i=\sum_{i=0}^4B_i=1 
$
for LPs of lengths $\ell = 5, 15, 25, 35, 45, 55, 65,75,85,95,105, 115$.
In Table~\ref{table:SOS}, we summarize the all-odd solutions to~(\ref{SOS_Diophantine_equations}) for odd $\ell$ such that $\ell=5m$ and $m \in \{1,\ldots,23\}$. In the last column, we record the all-odd solutions (up to sign changes and permutations) of the five sums of squares Diophantine equations of the fourth column. Since the linear equations $\sum_{i=0}^4A_i=\sum_{i=0}^4B_i=1$ must also be satisfied, some of these all-odd solutions are ruled out, and this is indicated by boldface. Conjecture~\ref{conj:Ilias} has practical value as it can  be used to prune the search space when $
\ell \equiv 0$ (mod $5$) by only decompressing $m$-compressions which are compatible with the all-odd solutions in Table~\ref{table:SOS}. 
\begin{table}[ht!]
\centering
{\scriptsize
$$
\begin{array}{c|c|c|c|c}
     \ell & m & \begin{array}{c} {\rm PSD}_\ba(rm) \\ {\rm PSD}_\bb(rm) \\ \end{array} &  \begin{array}{c}
          \sum_{i=0}^4A_i^2  \\
           \sum_{i=0}^4B_i^2 \\ 
     \end{array} & \mbox{All-odd solutions} \\
     \hline
     5  &  1   & 5 + 1 \pm (-1)^{\left\lfloor\frac{r}{2}\right\rfloor}\displaystyle \sqrt{5}
      x & 4 \cdot 1 + 1 = 5 & [ 1,1,1,1,1] \\
     15  &  3   & 15 + 1 \pm (-1)^{\left\lfloor\frac{r}{2}\right\rfloor}\displaystyle \sqrt{5}  x & 4 \cdot 3 + 1 = 13 & [ 1,1,1,1,3] \\
     25  &  5   & 25 + 1 \pm (-1)^{\left\lfloor\frac{r}{2}\right\rfloor}\displaystyle \sqrt{5}  x & 4 \cdot 5 + 1 = 21 & [ 1,1,1,3,3 ] \\
     35 &   7   & 35 + 1 \pm (-1)^{\left\lfloor\frac{r}{2}\right\rfloor}\displaystyle \sqrt{5}  x & 4 \cdot 7 + 1 = 29 & [1,1,1,1,5],\, [1,1,3,3,3]  \\
     45 &   9   & 45 + 1 \pm (-1)^{\left\lfloor\frac{r}{2}\right\rfloor}\displaystyle \sqrt{5}  x & 4 \cdot 9 + 1 = 37 & [1,1,1,3,5],\, [1,3,3,3,3]  \\
     55 &   11  & 55 + 1 \pm (-1)^{\left\lfloor\frac{r}{2}\right\rfloor}\displaystyle \sqrt{5}  x & 4 \cdot 11 + 1 = 45 & [1,1,3,3,5],\, \mathbf{[3,3,3,3,3]} \\
     65 &   13  & 65 + 1 \pm (-1)^{\left\lfloor\frac{r}{2}\right\rfloor}\displaystyle \sqrt{5}  x & 4 \cdot 13 + 1 = 53 & \mathbf{[1,1,1,1,7]},\, [1,1,1,5,5],\, [1,3,3,3,5] \\
     75 &   15  & 75 + 1 \pm (-1)^{\left\lfloor\frac{r}{2}\right\rfloor}\displaystyle \sqrt{5}  x & 4 \cdot 15 + 1 = 61 & [1,1,1,3,7],\,
     [1,1,3,5,5],\, [3,3,3,3,5] \\
     85 &   17  & 85 + 1 \pm (-1)^{\left\lfloor\frac{r}{2}\right\rfloor}\displaystyle \sqrt{5} x & 4 \cdot 17 + 1 = 69 & 
     [1,1,3,3,7],\, [1,3,3,5,5] \\
     95 &   19  & 95 + 1 \pm (-1)^{\left\lfloor\frac{r}{2}\right\rfloor}\displaystyle \sqrt{5}  x & 4 \cdot 19 + 1 = 77 & 
     [1, 1, 1, 5, 7],\, \mathbf{[1, 1, 5, 5, 5]},\, [1, 3, 3, 3, 7],\, [3, 3, 3, 5, 5] \\
     105 &   21  & 105 + 1 \pm (-1)^{\left\lfloor\frac{r}{2}\right\rfloor}\displaystyle \sqrt{5}  x & 4 \cdot 21 + 1 = 85 & \mathbf{[1, 1, 1, 1, 9]},\, [1, 1, 3, 5, 7],\, [1, 3, 5, 5, 5],\, [3, 3, 3, 3, 7] \\
     115 &   23  & 115 + 1 \pm (-1)^{\left\lfloor\frac{r}{2}\right\rfloor}\displaystyle \sqrt{5}  x & 4 \cdot 23 + 1 = 93 & 
     \mathbf{[1, 1, 1, 3, 9]},\, [1, 3, 3, 5, 7],\, [3, 3, 5, 5, 5]
\end{array} $$
\caption{All-odd solutions to equations~(\ref{SOS_Diophantine_equations})}
\label{table:SOS}
}
\end{table}

\section{Finding LPs of lengths $\ell = 85$ and  $\ell=87$}\label{sec:new85}
The orbit of a vector $\bv$ under circulant shifts
and decimations is called the {\em decimation class} of $\bv$~\cite{turner2020counting}. Since  $\ZZtimes$ acts on LPs, the search space for  LPs  is drastically reduced by searching only
across decimation class representatives~\cite{FGS:2001}. 
In fact, to classify $\{0,1\}$ LPs up to equivalence, in light of Lemma~\ref{lem:LPone}, Fletcher et al.~\cite{FGS:2001} exhaustively generated a set of all decimation class representatives 
$ \bv \in  \{0,1\}^{\ell}$  for odd $\ell \leq 47$ such that $\sum_{i=0}^{\ell-1}b_i=(\ell+1)/2$.
By Lemma~\ref{lem:LPone}, they then deleted the class representatives $\bv$ such that ${\rm PSD}_{\bv}(j)> (\ell+1)/2$ for some $j \in \zl\backslash\{0\}$. Then, among all the decimations  of a given vector $\bv$, they selected the
decimations $d_k(\bv)$ with $|{\rm PSD}_{d_k(\bv)}(1)-(\ell+1)/4|$ of greatest value to be the representatives
of  $\bv$'s decimation class, where  the decimation class of $\bv$ is allowed to have multiple different decimation class representatives.
In the second step, they sorted the list of class representative $\bv$s according to 
$${\rm PSD}_{\bv}(1)-\frac{\ell+1}{4}.$$
In the last step, they located  pairs of class representatives $\bu,\bv$ that satisfied
\begin{equation}\label{eqn:psdcrude}
{\rm PSD}_{\bu}(1)-\frac{\ell+1}{4}+{\rm PSD}_{\bv}(1)-\frac{\ell+1}{4}=0. 
\end{equation}
 To confirm that the LP property was satisfied by all the resulting pairs, they also deleted 
  pairs $(\bu,\bv)$ satisfying equation~(\ref{eqn:psdcrude}) that did not satisfy the equations involving the PSDs of $\bu$ and $\bv$ as in Lemma~\ref{lem:LPone}. In the end,  Fletcher et al.~\cite{FGS:2001} classified all LPs of length less than or equal to $47$ up to equivalence using this method.
A similar method that additionally exploited simultaneous  decompressions of candidate compressed vectors was implemented 
in~\cite{TKBG:DCC:2021}. This method found an LP  of length $77$ for the first time.

For  lengths $\ell=55$ and $77$,
Table~\ref{table:decsol} reports the decimation class counts for the number of decimation classes of length $\ell$, $\{0,1\}$ vectors with $(\ell+1)/2$ ones reported in~\cite{turner2020counting}, and the solution times for the simultaneous decompression-based searches reported in~\cite{TKBG:DCC:2021}.
\begin{table}[ht!]
\centering
    \begin{tabular}{c|c|c|c}
       Length  & Number of   & Exhaustive Search & One Solution \\
         $\ell$ &   Decimation Classes  &  (CPU Hours)& (CPU Hours)\\
        \hline 
        $55$ & $1.738341231644e+12$ & $101{,}542.4$ & --- \\
        $77$ & $2.945564382817e+18$ & --- &  $182{,}280$ \\
        \hline 
    \end{tabular} \caption{Decimation class counts in~\cite{turner2020counting} and solution times for  simultaneous decompression-based searches in~\cite{TKBG:DCC:2021}}
    \label{table:decsol}
    \end{table}

Assuming that the number of compressed pairs and the time required for decompression scales linearly on average with the number of decimation classes, we expect the Table~\ref{table:expected} time requirements for an exhaustive search or partial search to first solution if the method 
in~\cite{TKBG:DCC:2021} is used for lengths $\ell=85$ and $\ell=87$.
\begin{table}[ht!]
\centering
    \begin{tabular}{c|c|c|c}
       Length  & Number of   & Exhaustive Search & One Solution \\
         $\ell$ &   Decimation Classes  &  (CPU Hours)& (CPU Hours)\\
        \hline 
        $85$ & $6.100692175209e+20$ & $3.563621e+13$ & $37{,}752{,}839$ \\
        $87$ & $2.693812140345e+21$ & $1.573547e+14$ &  $166{,}700{,}847$ \\
        \hline 
    \end{tabular}
     \caption{Expected time needed  for an exhaustive search and partial search to first solution using the method 
in~\cite{TKBG:DCC:2021}}
\label{table:expected}
    \end{table}

These time estimates are optimistic as experiments suggest that the growth rate is superlinear in the number of decimation classes.  This growing magnitude of complexity precluded Turner~\etal~\cite{TKBG:DCC:2021}  from investigating larger open LP problems.  Partial searches for LPs using the compression method for $\ell>77$ only remain viable by reducing the number of complementary compressed vectors, or rather selecting the complementary compressions with higher likelihood of producing an LP. 
Alternatively, by Corollary~\ref{cor:GuhatGu} we can search for an LP $(\bu,\bv)$ that
has multiplier groups $\widehat{G}_{\bu}$ and $\widehat{G}_{\bv}$ both containing some non trivial elements of $\mathbb{Z}_n^{\times}$. This would reduce the search space drastically and can bring $\ell>77$ cases within computational reach if such an LP exists. For $\ell=3m$, $m \in \{39,43,49\}$, Kotsireas and Koutschan~\cite{KK:JCD:2021} used
 a method that exploits such a property along with a method that reduces the possible values of ${\rm PSD}_{\ba}(m)$ and keeps only the vectors whose ${\rm PSD}_{\ba}(m)$ are within the reduced set, 
and found the first examples of LPs of lengths $117$, $129$, and $147$. Kotsireas and Koutschan~\cite{KK:JCD:2021} also found an LP of length $133$ by solely restricting the search to LPs $(\bu,\bv)$ whose multiplier groups $\widehat{G}_{\bu}$ and $\widehat{G}_{\bv}$ both contained the group $G = \{1,11, 121\}<\mathbb{Z}_{133}^{\times}$.
\subsection{LPs of length $\ell = 85$}
We used the method in Section 4.1.1~in~\cite{KK:JCD:2021} with $H_1=\{ 1, 69 \}$ and Corollary~\ref{cor:GuhatGu} to restrict the search space for LPs of length $\ell = 85$ by assuming that the sought after LP $(\bu,\bv)$ has multiplier groups  $\widehat{G}_{\bu}$  and $\widehat{G}_{\bv}$
satisfying $H_1=\{ 1, 69 \}\leq\widehat{G}_{\bu}$ and
 $H_1=\{ 1, 69 \}\leq \widehat{G}_{\bv}$.
This method found the first known examples of  LPs of length $\ell = 85$. This had been the smallest previously unknown length case for LPs.  The subgroup $H_1= \{ 1, 69\}$ of $\mathbb{Z}_{85}^\times$  acts on $\mathbb{Z}_{85}$ and yields $16$ orbits of size $1$ and $34$ orbits of size $2$. We  searched for an LP of length $\ell=85$ which  could be obtained by combining the orbits of the subgroup  $H_1$. This restriction has the benefit of reducing the search space provided that such an LP exists.
  We chose $12$ orbits of size $1$ and $15$ orbits of size $2$ to  make  blocks of size $12 \cdot 1 + 15 \cdot 2 = 42$. Here, each block of size $42$ consists of positions of $-1$s determining the vectors $\bu$ and $\bv$. Therefore, the size of the search space was $\binom{16}{12} \cdot \binom{34}{15} = 1820 \cdot 1{,}855{,}967{,}520 = 3{,}377{,}860{,}886{,}400$. This search was not exhaustive and  was interrupted after traversing 2.6\% of the
search space. This took about 100 hours of CPU time. This search was done on 64 compute nodes, each with two 8-core Intel Haswell CPUs (Xeon
E5-2630v3, 2.4Ghz) and 128 GB RAM.
  For this search, we implemented the same PSD test in Section 4 of~\cite{KK:JCD:2021} without checking if the PSD value at $\ell/3$ is from a finite list of candidates, and used the same method as in Section 4 of~\cite{KK:JCD:2021} to identify LPs among vectors that passed the PSD test. However, this search did not take advantage of Conjecture~\ref{conj:Ilias} by only keeping the vectors $\ba$ that satisfied ${\rm PAF}_{{\mathcal A}_{17}}(0)=69$ for the $17$ compression  ${\mathcal A}_{17}$ of $\ba$ as Conjecture~\ref{conj:Ilias} had not yet been formulated at the time of the search. 
  Our search yielded $4$ equivalent ($1$ non-equivalent) LPs of length $\ell =85$ made out of $6$ different vectors. Their lexicographic rank encodings as subsets of size $12$  out of $16$ and $15$ out of $34$ are 
  \begin{align*}
  &\big (\{12, 1321116338\},\, \{42, 1275934280\}\big ), \\
  &\big (\{12, 1843909851\}, \, \{42, 606586783\}\big ), \\
  &\big (\{42, 1275934280\}, \, \{9, 1555522731\}\big), \\
  &\big (\{42, 606586783\},\, \{9, 788215097\}\big).
  \end{align*}
See~\cite{Kreher_Stinson} for lexicographic ranking and unranking algorithms for subsets of size~$k$. For the first LP $(\ba,\bb)$ of length $\ell = 85$ shown above, its  lexicographic rank encoding $\big(\{12, 13211163$ $38\},$  $\{42, 1275934280\}\big)$ is decoded as follows.
\begin{itemize}
    \item In the space of $\binom{16}{12} = 1820$ $12$-subsets of $\{1,\ldots,16\}$, decode 
    \begin{align*}
        & 12 \mbox{ as }  \ba_\mathrm{ones} = \{1, 2, 3, 4, 5, 6, 7, 8, 9, 10, 14, 15\}, \\ 
        & 42 \mbox{ as }  \bb_\mathrm{ones} = \{1, 2, 3, 4, 5, 6, 7, 8, 10, 11, 14, 15\}. 
    \end{align*}
    \item In the space of $\binom{34}{15} = 1,855,967,520$ $15$-subsets of $\{1,\ldots,34\}$, decode
        \begin{align*}
        & 1321116338 \mbox{ as }  \ba_\mathrm{twos} = \{3, 4, 5, 7, 10, 11, 22, 24, 25, 27, 28, 29, 30, 31, 34\}, \\ 
        & 1275934280 \mbox{ as }  \bb_\mathrm{twos} = \{2, 8, 10, 11, 12, 15, 19, 21, 23, 25, 26, 28, 29, 33, 34\}. 
    \end{align*}
    \item Enumerate the $16$ orbits of size $1$ in increasing order as \\ 
    $O_1 = \{ \{5\}, \{10\}, \{15\}, \{20\}, \{25\}, \{30\}, \{35\}, \{40\}, \{45\}, \{50\}, \{55\}, \{60\}, \{65\},$ \\
    $\{70\}, \{75\}, \{80\} \}.$
    \item Enumerate the $34$ orbits of size $2$ in increasing order of their smallest element as \\ 
    $O_2 = \{ \{1, 69\}, \{2, 53\}, \{3, 37\}, 
    \{4, 21\}, \{6, 74\}, \{7, 58\}, \{8, 42\}, \{9, 26\},
    \{11, 79\}, $ \\ $ \{12, 63\},  
    \{13, 47\}, \{14, 31\}, \{16, 84\}, \{17, 68\}, \{18, 52\}, \{19, 36], \{22, 73], \{23, 57\},$ \\ $  \{24, 41\}, \{27, 78\}, 
      \{28, 62\}, \{29, 46\},\{32, 83\},\{33, 67\}, \{34, 51\}, \{38, 72\}, \{39, 56\}, $
      \\  $   \{43, 77\}, \{44, 61\}, \{48, 82\},    \{49, 66\},\{54, 71\},\{59, 76\},\{64, 81\} \}.$ 
    \item Make the block of size $42$ of the indices of the positions of the $-1$ elements in $\ba$, by combining $12$ elements of $O_1$ whose indices are given by $\ba_\mathrm{ones}$ and $15$ elements of $O_2$ whose indices are given by $\ba_\mathrm{twos}$. This yields the following $\ba$-block of size $42$: \\
    $\{5, 10, 15, 20, 25, 30, 35, 40, 45, 50, 70, 75, 3,
    37, 4, 21, 6, 74, 8, 42, 12, 63, 13, 47, 29, 46, $ \\  $33, 67, 34, 51, 39, 56, 43,
    77, 44, 61, 48, 82, 49, 66, 64, 81\}$. 
    \item Make the block of size $42$ of the indices of the positions of the $-1$ elements in~$\bb$, by combining $12$ elements of $O_1$ whose indices are given by $\bb_\mathrm{ones}$ and $15$ elements of $O_2$ whose indices are given by $\bb_\mathrm{twos}$.
    This yields the following $\bb$-block of size $42$: \\
    $\{5, 10, 15, 20, 25, 30, 35, 40, 50, 55, 70, 75, 2,
    53, 9, 26, 12, 63, 13, 47, 14, 
    31, 18, 52, 24,$ \\
    $41, 28, 62, 32, 83, 34, 51, 38, 72, 43, 77, 44, 61, 59, 76, 64, 81\}$. 
    \item The above $\ba$-block and $\bb$-block (both of size 42) 
    yield an LP for $\ell = 85$.
\end{itemize}
 For the first LP $(\ba,\bb)$ of length $\ell = 85$ shown above,
 the $17$-compressions are
$$
\mathcal{A}_{17} = [1, 3, 3, 1, -7], \, \,\,\mathcal{B}_{17} = [3, 1, 1, 3, -7].
$$
These possess the properties required by Proposition~\ref{proposition_1}. That is, 
\[{\rm PAF}_{\mathcal{A}_{17}}(0) = {\rm PAF}_{\mathcal{B}_{17}}(0) = 4 \cdot 17 + 1 = 69.\] 
As a consequence, the coefficients of $\displaystyle \sqrt{5}$ in ${\rm PSD}_{\ba}(17r)$ and ${\rm PSD}_{\bb}(17r)$ cancel out for $r=1,\ldots,4$. Specifically,
$$
\begin{array}{l}
{\rm PSD}_{\ba}(17r)  =  85 +1 +(-1)^{\left\lfloor \frac{r}{2} \right\rfloor}\displaystyle \sqrt{5} \cdot 18   \\[2ex] 
{\rm PSD}_{\bb}(17r)  =  85 +1 -(-1)^{\left\lfloor \frac{r}{2} \right\rfloor}\displaystyle\ \sqrt{5} \cdot 18. 
\end{array}
$$
\begin{remark}
 We did not use Conjecture~\ref{conj:Ilias} or Corollary~\ref{cor:k1plusk2} in our search for an  LP of length $85$ because neither Conjecture~\ref{conj:Ilias} nor Corollary~\ref{cor:k1plusk2} had been formulated when we implemented our search.    However, the LP we found happened to satisfy 
 equations~(\ref{eqn:Ilias}) in Conjecture~\ref{conj:Ilias}.
\end{remark}
\subsection{LPs of length $\ell = 87$}

Next, we describe our search method for LPs $(\ba,\bb)$ of length $\ell=87$.
First we need the following theorem which is a special case of  Theorem 3 in~\cite{DK:DCC:2015}.
\begin{theorem}\label{thm:Dragomir}
 Let $(\ba,\bb)$ be an LP of length $\ell=dm$. Then the  $m$-compressed vectors  $(\mathcal{A}_{m},\mathcal{B}_{m})$ satisfy
  \begin{align*}
{\rm PAF}_{\mathcal{A}_{m}}(0)+{\rm PAF}_{\mathcal{B}_{m}}(0)&=2\ell-2\left(\frac{\ell}{d}-1\right)=2\ell+2-\frac{2\ell}{d},\label{eqn:PAFcompconstraint1}\\
{\rm PAF}_{\mathcal{A}_{m}}(j)+{\rm PAF}_{\mathcal{B}_{m}}(j)&=-\frac{2\ell}{d}, \quad \forall j\in \Z_d-\{0\}.
\end{align*}
\end{theorem}
By Theorem~\ref{thm:Dragomir},  the  $3$-compressed vectors  $(\mathcal{A}_{3},\mathcal{B}_{3})$ of an LP $(\ba,\bb)$ of length $\ell = 87$ must contain $14$ elements with absolute value equal to $3$ and $2 \cdot 29 - 14 = 44$ elements with absolute value equal to $1$. Experimental evidence from the study of other lengths 
which are divisible by $3$ indicates that 
$\mathcal{A}_{3}$ and $\mathcal{B}_{3}$ which have an equal number of elements with absolute value equal to $3$ are more likely to yield LPs. 
Hence, we computed approximately six thousand candidate $3$-compressions 
satisfying the following constraints. 
\begin{enumerate}
    \item The vectors $\mathcal{A}_{3}, \mathcal{B}_{3} \in \{-3,-1,+1,+3\}^{29}$ contain $14$ elements with absolute value equal to $3$ and  $44$ elements with absolute value equal to $1$.
    \item ${\rm \PAF}_{\mathcal{A}_{3}}(s) + {\rm \PAF}_{\mathcal{B}_{3}}(s)  = (-2) \cdot 3 = -6$ for $s=1,\ldots,28$.
    \item ${\rm \PSD}_{\mathcal{A}_{3}}(s) + {\rm \PSD}_{\mathcal{B}_{3}}(s)  = 2\ell+2= 2 \cdot 87 + 2 = 176$ for $s=1,\ldots,28$.
    \item $\sum_{i=0}^{28}A_i=\sum_{i=0}^{28}B_i=1$. 
    \item  Each of $\mathcal{A}_{3}$ and $\mathcal{B}_{3}$ contains $7$ elements with absolute value equal to $3$ and $22$ elements with absolute value equal to $1$.
\end{enumerate}
\begin{singlespace}
{\normalsize Among the above constraints, only  constraint 5 is not necessary for a length $87$ LP to exist. However, imposing constraint 5 was essential in greatly reducing the search space to a part where solutions are most likely to exist. 
Subsequently, we ran our {\tt C} $3$-uncompression code for approximately two thousand 
candidate $3$-compressions and discovered the following two LPs of order $\ell = 87$:}
\end{singlespace}
{\footnotesize
\begin{align*}
\ba &= [-1,-1,-1,-1,1,-1,-1,-1,1,1,-1,1,-1,1,-1,1,1,-1,-1,1,1,-1,-1,1,-1,1,-1,1,1, \\
&\qquad {-1},-1,1,1,1,-1,-1,-1,-1,-1,1,-1,-1,1,-1,-1,-1,1,1,1,1,-1,-1,-1,1,-1,1,1,\\
&\qquad 1,-1,1,1,-1,-1,-1,-1,1,1,1,1,1,1,-1,-1,1,1,1,-1,1,1,1,1,1,-1,1,-1,1,-1],\\
\bb &= [-1,-1,-1,-1,-1,-1,-1,-1,1,1,1,1,-1,1,1,-1,1,1,-1,1,-1,1,1,1,1,-1,-1,1,1,\\
&\qquad {-1},1,1,-1,1,-1,-1,1,-1,-1,1,1,1,-1,1,-1,1,1,-1,-1,1,-1,-1,1,1,1,1,1,-1,-1,\\
&\qquad {-1},1,1,1,-1,-1,1,-1,1,-1,-1,1,-1,1,-1,-1,-1,1,-1,-1,1,1,1,-1,1,-1,1,-1],\\
\ba &= [-1,-1,-1,-1,1,-1,-1,-1,1,1,1,1,1,1,-1,1,-1,-1,-1,1,1,1,-1,1,-1,-1,-1,1,\\
&\qquad 1,-1,-1,1,1,1,-1,-1,-1,1,-1,-1,1,-1,1,-1,-1,1,1,1,1,1,-1,-1,1,1,1,-1,1,1,-1,\\&\qquad 1,1,-1,-1,-1,-1,1,-1,1,1,-1,-1,-1,-1,1,1,1,-1,1,1,-1,1,-1,-1,1,1,1,-1],\\
\bb &= [-1,1,-1,1,1,-1,-1,1,1,1,1,1,1,-1,1,-1,1,1,1,-1,-1,1,-1,1,-1,1,-1,1,-1,\\
&\qquad {-1},-1,1,-1,-1,-1,-1,-1,-1,-1,1,-1,1,1,1,-1,1,-1,-1,1,-1,-1,1,1,1,-1,-1,\\
&\qquad 1,1,-1,-1,1,-1,1,-1,-1,1,-1,1,-1,1,-1,-1,1,-1,-1,1,-1,-1,1,1,1,1,1,1,1,1,-1].
\end{align*}
}
Both of the above LPs of length $\ell = 87$ $3$-compress to 
{\small
\begin{align*}
\mathcal{A}_{3} &= [-3, -1, 1, -1, 1, -3, -3, -1, 1, 1, 1, 1, -1, 1, -3, 1, 1, 1, -1, 3, 3, -1, -1, 1, -1, 1, -1, 3, 1],\\
\mathcal{B}_{3} &= [-3, -1, 1, -1, 1, -3, -3, 1, -1, 1, 1, 1, 1, -1, 3, -3, 1, 1, -1, -1, -1, 1, 1, 3, 1, 1, -1, 3, -1].
\end{align*}
}
There are seven $\pm 3$s in each of the two vectors $\mathcal{A}_{3}$ and $\mathcal{B}_{3}$, 
and this ``balanced" configuration  yields LPs of length $\ell = 87$.
Based on the analysis at the beginning of this section, we can make the following claim.
\begin{claim}
Without imposing constraint 5, a successful outcome  of our search for a length $87$ LP would not have been possible.
\end{claim}
\section{Searching for LPs of length $\ell=115$}\label{sec:search115}
By using Corollary~\ref{cor:GuhatGu}, two non-exhaustive searches for an LP $(\bu,\bv)$ of length $115$ with multiplier groups $\widehat{G}_{\bu} $ and $\widehat{G}_{\bv} $ both containing the subgroup $\{1,91\}$ were performed.
The first search was done before  Conjecture~\ref{conj:Ilias} was formulated, hence
it did not use it. For the first search, the combinations of $0$ cosets of size $1$ and $29$
cosets of size $2$ were considered, yielding a search space of size
   $\binom{4}{0}  \cdot \binom{55}{29} = 3{,}560{,}597{,}348{,}629{,}860$.
The computation was aborted after almost $4\%$ of the search space were
traversed. This took about $4,359$ days of CPU time, and the output files
took up $76$ GB of disk space. 
The second search was done with the implementation of Conjecture~\ref{conj:Ilias}. This
time, combinations of $4$ cosets of size $1$ and $27$ cosets of size $2$ were
considered, yielding a search space of size
   $\binom{4}{4}  \cdot \binom{55}{27} = 3{,}824{,}345{,}300{,}380{,}220.$
This computation was aborted after $3{,}436$ CPU days, when slightly more
than $10\%$ of the search space was traversed. The output files occupied
about $64$ GB of disk space. For both of the searches,  the same PSD test in Section 4 of~\cite{KK:JCD:2021} without checking if the PSD value at $\ell/3$ is from a finite list of candidates,
was implemented. Both searches used the same method as in Section 4 of~\cite{KK:JCD:2021} to identify LPs among vectors that passed the PSD test.
The comparison of the two partial searches shows the significant gain (both time-wise
and space-wise) that is obtained from implementing Conjecture~\ref{conj:Ilias}.
Neither the first, nor the second search yielded an LP of length $115$. Both of the searches were done on 64 compute nodes, each with two 8-core Intel Haswell CPUs (Xeon
E5-2630v3, 2.4Ghz) and 128 GB RAM.
\section{Conclusion and future research}\label{sec:conclusion}
Recently, LPs of 
length $77$ were found in~\cite{TKBG:DCC:2021}, and lengths $117, 129, 133, 147$ were found  in~\cite{KK:JCD:2021} for the first time. In this  paper, we find LPs of (the previously open) lengths $85, 87$. 
This reduces the list of integers less than $200$ for which the existence of LPs problem  remains open to the following ten values:
$$
115, 145,  159, 161, 169, 175, 177, 185, 187, 195.
$$

An LP $(\ba, \bb)$ as defined in this article corresponds to a difference family in 
$\mathbb{Z}_{\ell}$~\cite{Hollon}.
In~\cite{Iliasspecial:2019},  previously known theory  to search for
difference families in $\mathbb{Z}_{\ell}$ was generalized  to search for
difference families in finite abelian groups. A possible direction for future research is generalizing the 
theory in this paper to difference families in finite abelian groups, i.e., by the structure theorem of finite abelian groups, groups of the form
$\mathbb{Z}_{d_1} \times \mathbb{Z}_{d_2}\times\cdots \times \mathbb{Z}_{d_m} $, where $d_i$ divides $d_{i+1}$ for $i=1,\ldots,m-1$, and $d_m\equiv 0 \pmod 5$.
\section*{Acknowledgments}
  The views expressed in this article are those of the authors, and do not reflect the official policy or position of the United States Air Force, Department of Defense, or the U.S.~Government.  
  Christoph Koutschan was supported by the Austrian Science Fund (FWF): F5011-N15. The authors thank three referees for carefully reading the paper and providing  comments that significantly improved the paper.
  
  Data Deposition Information: No data sets have been used.
  
\bibliographystyle{plain}
\bibliography{Legendre_pairs}
\end{document}